\documentclass{mcom-l}

\usepackage[english]{babel}
\usepackage{graphics}
\usepackage{graphicx}
\usepackage{epsfig}
\usepackage{amssymb}
\usepackage{arydshln}
\usepackage{amsmath}
\usepackage{amsfonts}
\usepackage{float}

   \newtheorem{theorem}{Theorem}[section]
   \newtheorem{prop}[theorem]{Proposition}

   \theoremstyle{definition}
     \newtheorem{definition}[theorem]{Definition}
    \newtheorem{example}[theorem]{Example}
    
    \theoremstyle{remark}
   \newtheorem{remark}[theorem]{Remark}

  \numberwithin{equation}{section}

\newcommand{\B}[1]{\mbox{\boldmath $#1$}}
\newcommand{\be}{\begin{equation}}
\newcommand{\ee}{\end{equation}}
\newcommand{\ba}{\begin{array}}
\newcommand{\ea}{\end{array}}

\newcommand{\al}{{\alpha}}
\newcommand{\bt}{{\beta}}

\newcommand{\g}{{\gamma}}

\newcommand{\La}{\Lambda}

\newcommand{\Om}{\Omega}

\newcommand{\sg}{\sigma}
\newcommand{\tl}{\tilde}

\begin{document}

\title{Implicit QR for companion-like  pencils}

\author{P. Boito}
\address{XLIM--DMI, UMR CNRS 7252, Facult\'e des Sciences et 
Techniques,
123 av. A. Thomas, 87060 Limoges, France}
\email{paola.boito@unilim.fr}

\author{Y. Eidelman}
\address{School of Mathematical Sciences, Raymond and Beverly
Sackler Faculty of Exact Sciences, Tel-Aviv University, Ramat-Aviv,
69978, Israel}
\email{eideyu@post.tau.ac.il}

\author{L. Gemignani}
\address{Dipartimento di Informatica, Universit\`a di Pisa,
Largo Bruno Pontecorvo 3, 56127 Pisa, Italy}
\email{l.gemignani@di.unipi.it}
\thanks{This work was partially supported by MIUR, grant number
20083KLJEZ.}

\subjclass[2010]{65F15, 65H17}

\date{Jan 29, 2014}


\begin{abstract}
A fast implicit QR algorithm for eigenvalue computation of
low   rank   corrections  of unitary  matrices is  adjusted to work with
matrix pencils arising from  polynomial zerofinding problems .  The modified  
QZ  algorithm computes the
generalized   eigenvalues of certain  $N\times N$
 rank structured matrix pencils
 using $O(N^2)$ flops and $O(N)$ memory storage.
Numerical experiments  and comparisons confirm the 
effectiveness and the
stability of the proposed  method.
\end{abstract}

%

\maketitle

\section{Introduction}
Computing the roots of a univariate polynomial is a fundamental problem that arises in many applications. 
One way of numerically computing the roots of a polynomial is to form its companion matrix (pencil) and 
compute the (generalized) eigenvalues. 

In a paper on polynomial root--finding \cite{JV}, J\'onsson and Vavasis present a comparative analysis 
of the accuracy of different matrix algorithms. The conclusion is that computing the roots of a polynomial 
by first forming the associated companion pencil $A-\lambda B$ and then solving $A\B x=\lambda B\B x$ using the 
QZ algorithm provides better backward error bounds than computing the eigenvalues of the associated 
companion matrix by means of the QR algorithm.   The analysis does not take into account the possible use of 
balancing which would have the effect  of allineating  the accuracy of the two approaches \cite{TT,DR}. However, 
it is worth noting that this use is not allowed if the customary QR and QZ algorithms are modified 
at the aim of reducing their complexity by one order of magnitude by 
exploiting  the rank structure of the initial companion matrix (pencil). Thus, the use of the QZ algorithm 
applied to the companion pencil achieves  the potential best score in terms of  accuracy and efficiency 
among matrix methods for  polynomial root--finding.

Another interesting  appearance of the pencil approach is in the design of 
numerical routines for 
computing the determinant hence the zeros of polynomial matrices \cite{POLY}. 
A classical technique is the interpolation of the determinant of 
the polynomial matrix $A(\lambda)$  at the roots of unity via FFT. This gives a representation of 
$p(\lambda):=
\det A(\lambda)$ in the basis of Lagrange polynomials generated from the set of nodes. 
The main stumbling step in this process is the presence of undesirable infinite zeros or, equivalently, 
zero leading coefficients of $p(\lambda)$.   This makes the companion matrix 
approach fully unusable. On the contrary, a companion pencil can still be formed 
whose generalized eigenvalues are the 
finite roots of the polynomial.  If the pencil is $A-\lambda B$ directly constructed 
from the coefficients of the 
polynomial expressed in the Lagrange basis then $A$ is a unitary plus rank--one matrix and $B$ is a rank--one 
perturbation of the identity matrix. By a preliminary reduction using the algorithm in \cite{egg_qr}
the matrices $A$ and $B$ can be transformed in upper Hessenberg and triangular form, respectively. 
  
This paper undertakes the task of developing  an efficient 
variant of the QZ algorithm applied to  a generalized companion pair $(A, B)$,
where $A$ is upper Hessenberg, $B$ is upper triangular and  $A$ and $B$ are rank--one 
modifications of unitary matrices.  These assumptions  imply  
suitable rank structures in both $A$ and $B$. 
We show that  each matrix  pair generated by the
QZ process
is also  rank-structured.  By exploiting  this property 
a  novel,    fast adaptation   of the  QZ
eigenvalue algorithm    is obtained
 which  requires  $O(n)$   flops  per iteration
and $O(n)$  memory space only.
Numerical experiments confirm that the algorithm
is stable.

The paper is organized as follows. In Sect.~2 we introduce the computational
problem and
describe the  matrix structures involved.   Sect.~3 and  Sect.~4  deal with
basic issues and concepts concerning the  condensed representation  and the
invariance of these structures under the QZ process, whereas the main algorithm is presented in 
Sect.~5. Sect.~6 gives a backward error analysis for the QZ method applied to the companion pencil. In Sect.~7 we present an implementation of our fast variant of the QZ algorithm   together
with the results of extensive numerical
experiments. Finally, the conclusion and a discussion
are the subjects of Sect.~8.

\section{The Problem Statement}
\setcounter{equation}{0}
Companion pencils and generalized companion pencils expressed in the 
Lagrange basis  at the roots of unity are specific instances of the 
following general class. 

\begin{definition}\label{pn}
The matrix pair $(A,B)$,  $A, B\in \mathbb C^{N\times N}$, belongs to the
class 
$\mathcal P_N \subset \mathcal C^{N\times N} \times \mathcal C^{N\times N}$   
of generalized companion pencils 
iff:
\begin{enumerate}
\item $A\in \mathbb C^{N\times N}$ is upper Hessenberg;
\item $B\in \mathbb C^{N\times N}$ is upper triangular;
\item There exist two vectors $\B z\in \mathbb C^N$ and $\B w\in \mathbb C^N$ 
and a unitary matrix $V\in \mathbb C^{N\times N}$  such that
\begin{equation}\label{a}
 A=V-\B z \B w^*; 
\end{equation}
\item There exist two vectors $\B p\in\mathbb C^N$ and $\B q\in\mathbb C^N$ and
a unitary matrix $U\in \mathbb C^{N\times N}$  such that
\begin{equation}\label{b}
 B=U-\B p \B q^*.
\end{equation}
\end{enumerate}
\end{definition}

In order to characterize the individual properties of the matrices $A$ and $B$ 
we give some additional definitions.

\begin{definition}\label{bn}
We denote  by  ${\mathcal T}_{N}$ the class of  upper triangular matrices 
$B\in \mathbb C^{N\times N}$ which are rank one perturbations of unitary
matrices, i.e., such that \eqref{b} holds for a suitable unitary matrix $U$ and
vectors $\B p, \B q$.
\end{definition}

Since $B$ is upper triangular the strictly lower triangular part of the unitary
matrix $U$ in \eqref{b} coincides with the corresponding part of the rank one 
matrix $\B p \B q^*$, i.e.,
\be\label{semop23}
U(i,j)= p(i) q^*(j),\quad 1\le j<i\le N,
\end{equation}
where $\{p(i)\}_{i=1,\ldots,N}$ and $\{q(j)\}_{j=1,\ldots,N}$ are the entries 
of $\B p$ and $\B q$, respectively.

\begin{definition}\label{un}
We denote by ${\mathcal U}_{N}$ the class of unitary matrices 
$U\in \mathbb C^{N\times N}$ which satisfy the condition \eqref{semop23}, i.e.,
for which there exist vectors $\B p, \B q$ such that
the matrix $B=U-\B p \B q^*$ is an upper triangular matrix.
\end{definition}

Observe that we have
\[
U\in {\mathcal U}_{N} \Rightarrow
{\rm rank} U(k+1\colon N,1\colon k)\le 1, \quad k=1,\dots,N-1.
\]
From the nullity theorem \cite{FM}, see also \cite[p.142]{EGH1} it follows 
that the same  property also holds in the strictly upper triangular part, 
namely,
\begin{equation}\label{setp26}
U\in {\mathcal U}_{N} \Rightarrow
{\rm rank} U(1\colon k, k+1 \colon N)\le 1, \quad k=1,\dots,N-1.
\end{equation}

\begin{definition}\label{hn}
We denote  by  ${\mathcal H}_{N}$ the class of  upper Hessenberg  matrices 
$A\in \mathbb C^{N\times N}$ which are rank one perturbations of unitary
matrices, i.e., such that \eqref{a} holds for a suitable unitary matrix $V$ and
vectors $\B z, \B w$.
\end{definition}

\begin{definition}\label{qn}
We denote by ${\mathcal V}_{N}$ the class of unitary matrices 
$V\in \mathbb C^{N\times N}$ for which there exist vectors $\B z,\B w$ such 
that the matrix $A=V-\B z \B w^*$ is an upper Hessenberg matrix.
\end{definition}

We find that 
\[
V\in {\mathcal V}_{N} \Rightarrow
{\rm rank} V(k+2\colon N,1\colon k)\le 1, \quad k=1,\dots,N-2.
\]
Again from the nullity theorem  it follows that a similar  property also
holds in the upper triangular part, namely,
\begin{equation}\label{setp27}
V\in {\mathcal V}_{N} \Rightarrow
{\rm rank} V(1\colon k, k \colon N)\le 2, \quad k=1,\dots,N.
\end{equation}

The QZ algorithm  is  the customary method for solving generalized eigenvalue 
problems numerically  by means of unitary  transformations
(see e.g. \cite{GVL} and \cite{WD}). Recall that the Hessenberg/triangular 
form is preserved under the QZ iteration; an easy computation then yields
\begin{equation}\label{qzss}
(A,B)\in \mathcal P_N, \ (A,B)\overset{{\rm QZ \ step}} 
\rightarrow (A_1,B_1) \Rightarrow (A_1,B_1)\in \mathcal P_N.
\end{equation}
Indeed if $Q$ and $Z$ are unitary then from (\ref{a}) and (\ref{b}) it follows
that the matrices $A_1=Q^*AZ$ and $B_1=Q^*BZ$ satisfy the relations
$$
A_1=V_1-\B z_1\B w_1^*,\quad B_1=U_1-\B p_1\B q_1^*
$$
with the unitary matrices $V_1=Q^*VZ,\;U_1=Q^*UZ$ and the vectors
$\B z_1=Q^*z,\;\B w_1=Z^*w,\;\B p_1=Q^*p,\;\B q_1=Z^*q$. Moreover one can 
choose the unitary matrices $Q$ and $Z$ such that the matrix $A_1$ is upper
Hessenberg and the matrix $B_1$ is upper triangular.
Thus, one  can in principle think of  designing a structured QZ iteration 
that, given in input a condensed representation of the matrix pencil  
$(A,B)\in \mathcal P_N$, returns as output a condensed representation of  
$(A_1,B_1)\in \mathcal P_N$ generated by one step of the classical QZ algorithm
applied to $(A,B)$. In the next sections we first  introduce an  eligible  
representation  of a  rank-structured matrix pencil $(A,B)\in \mathcal P_N$ and
then  discuss the modification of this representation under the QZ process.

\section{Quasiseparable Representations}
\setcounter{equation}{0}
In this  section we present the properties of  quasiseparable 
representations of rank--structured matrices \cite{EGfirst},
\cite[Chapters 4,5]{EGH1}.
First  we recall some general  results and  definitions. Subsequently, 
we describe their  adaptations 
for the representation of the matrices involved in the structured QZ 
iteration applied to an input 
matrix pencil $(A,B)\in \mathcal P_N$.

\subsection{{\bf Preliminaries}}

A matrix $M=\{M_{ij}\}_{i,j=1}^N$ is {\em $(r^L,r^U)$-quasiseparable}, with 
$r^L, r^U$ positive integers, if, using MATLAB\footnote{MATLAB is a registered 
trademark of The Mathworks, Inc..} notation, 
\begin{gather*}
\max_{1\leq k\leq N-1}{\rm rank}(M(k+1:N,1:k))\leq r^L,\\
\max_{1\leq k\leq N-1}{\rm rank}(M(1:k,k+1:N))\leq r^U.
\end{gather*}
Roughly speaking, this means that every submatrix extracted from the lower 
triangular part of $M$ has rank at most $r^L$, and every submatrix extracted 
from the upper triangular part of $M$ has rank at most $r^U$. Under this 
hypothesis, $M$ can be represented using $\mathcal{O}(((r^L)^2+(r^U)^2)N)$ parameters. 
In this subsection we present such a representation.

The quasiseparable representation of a rank--structured matrix consists of a 
set of vectors and matrices used to generate its  entries. For the sake of 
notational simplicity, generating  matrices  and vectors are  denoted by a 
lower-case  letter.
 
In this representation, the entries of $M$ take the form
 
 \begin{equation}\label{qs1}
M_{ij}=\left\{\begin{array}{ll}
p(i)a_{ij}^{>}q(j),&1\le j<i\le N,\\
d(i),&1\le i=j\le N,\\
g(i)b_{ij}^{<}h(j),&1\le i<j\le N\end{array} \right.
\end{equation}
 where:
 \begin{itemize}
 \item[-] $p(2),\ldots,p(N)$ are row vectors of length $r^L$,
  $q(1),\ldots,q(N-1)$ are column vectors of length $r^L$, and
  $a(2),\ldots,a(N-1)$ are matrices of size $r^L\times r^L$; these are called 
{\em lower quasiseparable generators} of order $r^L$;
 \item[-] $d(1),\ldots,d(N)$ are numbers (the diagonal entries),
 \item[-] $g(2),\ldots,g(N)$ are row vectors of length $r^U$,
  $h(1),\ldots,h(N-1)$ are column vectors of length $r^U$, and
  $b(2),\ldots,b(N-1)$ are matrices of size $r^U\times r^U$; these are called 
{\em upper quasiseparable generators} of order $r^U$;
\item[-]
 the matrices $a_{ij}^{>}$ and $b_{ij}^{<}$ are defined as
 \[
\left\{\begin{array}{ll}a_{ij}^{>}=a(i-1)\cdots
a(j+1) \  {\rm for}\  i>j+1;\\ a_{j+1,j}^{>}=1
\end{array}\right.
\]
and
\[
\left\{\begin{array}{ll}
b_{ij}^{<}=b(i+1)\cdots
b(j-1) \ {\rm  for} \  j>i+1; \\
b_{i,i+1}^{<}=1.
\end{array}\right.
\]
 \end{itemize}

From (\ref{setp26}) it follows that any matrix from the class 
${\mathcal U}_{N}$ has upper quasiseparable generators with orders equal one.

 The quasiseparable representation can be generalized to the case where $M$ is 
a block matrix, and to the case where the generators do not all have the same 
size, provided that their product is well defined. Each block $M_{ij}$ of size 
$m_i\times n_j$ is represented as in \eqref{qs1}, except that the sizes of the 
generators now depend on $m_i$ and $n_j$, and possibly on the index of $a$ and 
$b$. More precisely:
 \begin{itemize}
 \item[-]
 $p(i), q(j), a(k)$ are matrices of sizes $m_i\times
r^L_{i-1},\;r^L_j\times n_j,\; r^L_k\times r^L_{k-1}$,
respectively; 
\item[-]
$d(i)\;(i=1,\dots,N)$ are $m_i\times n_i$ matrices,
\item[-]
$g(i),h(j),b(k)$ are matrices of sizes
$m_i\times r^U_i,\;r^U_{j-1}\times n_j,\; r^U_{k-1}\times r^U_k$,
respectively.
 \end{itemize}
 The numbers $r^L_k,r^U_k\;(k=1,\dots,N-1)$ are called the {\it orders}
of these generators.
 
It is worth noting that lower and upper quasiseparable  generators of a matrix 
are not uniquely defined. A set of generators with minimal orders can be 
determined according to the ranks of maximal submatrices  located in  the lower
and upper triangular parts of  the matrix.

One advantage of the block representation for the purposes of the present paper
consists in the fact that $N\times N$ upper Hessenberg matrices can be treated 
as $(N+1\times N+1)$ block upper triangular ones by choosing blocks of sizes 
\be\label{aprmn18s}
m_1=\dots=m_N=1,\;m_{N+1}=0,\quad
n_1=0,\;n_2=\dots=n_{N+1}=1.
\end{equation}
Such a treatment allows also to consider quasiseparable representations which 
include the main diagonals of matrices. Assume that $C$ is an $N\times N$ 
scalar
matrix with the entries in the upper triangular part represented in the form
\be\label{appr18}
C(i,j)=g(i)b^<_{i-1,j}h(j),\quad 1\le i\le j\le N
\end{equation}
with matrices $g(i),h(i)\;(i=1,\dots,N),\;b(k)\;(k=1,\dots,N-1)$ of sizes
$1\times r_i,r_i\times1,r_k\times r_{k+1}$. The elements $g(i),h(i)\;
(i=1,\dots,N),\;b(k)\;(k=1,\dots,N-1)$ are called {\it upper triangular
generators} of the matrix $C$ with orders $r_k\;(k=1,\dots,N)$. 
From (\ref{setp27}) it follows that any matrix from the class 
${\mathcal V}_{N}$ has upper triangular generators
with orders not greater than two. 
If we treat a matrix $C$ as a block one with entries of sizes (\ref{aprmn18s}) 
we conclude
that the elements $g(i)\;(i=1,\dots,N),\;h(j-1)\;(j=2,\dots,N+1),\;b(k-1)\;
(k=2,\dots,N)$ are upper quasiseparable generators of $C$.

Matrix operations involving zero-dimensional arrays (empty matrices) are 
defined according to the rules used in MATLAB and described in \cite{deBoor}.
In particular, the  product of
a  $m\times 0$ matrix by a $0\times m$ matrix is a $m\times m$ matrix with all 
entries equal to 0.
Empty matrices may be used in assignment statements as a convenient way to add 
and/or delete rows or columns of matrices.

\subsection{{\bf Representations of matrix pairs from the class
${\mathcal P}_{N}$}}\label{psec}

Let $(A,B)$ be a matrix pair from the class ${\mathcal P}_{N}$.
The corresponding matrix $A$ from the class ${\mathcal H}_{N}$ is completely
defined by the following parameters:
\begin{enumerate}
\item  the subdiagonal entries $\sg^A_k\;(k=1,\dots,N-1)$ of the matrix $A$;  
\item the upper triangular generators $g_V(i),h_V(i)\;(i=1,\dots,N),\;
b_V(k)\;(k=1,\dots,N-1)$ of the corresponding unitary matrix $V$ from the class
${\mathcal V}_N$;
\item  the vectors of perturbation $\B z={\rm col}(z(i))_{i=1}^N,\;
\B w={\rm col}(w(i))_{i=1}^N$.
\end{enumerate}
From  \eqref{setp27} it follows that 
the matrix $V\in \mathcal V_N$ has upper triangular generators with orders
not greater than two.

The corresponding matrix $B$ from the class ${\mathcal T}_{N}$ is completely
defined by the following parameters:
\begin{enumerate}
\item  the diagonal entries $d_B(k)\;(k=1,\dots,N)$ of the matrix $B$;  
\item the upper quasiseparable generators $g_U(i)\;(i=1,\dots,N-1),\;
h_U(j)\;(j=2,\dots,N),\;b_U(k)\;(k=2,\dots,N-1)$ of the corresponding unitary 
matrix $U$ from the class ${\mathcal U}_N$;
\item  the vectors of perturbation $\B p={\rm col}(p(i))_{i=1}^N,\;
\B q={\rm col}(q(i))_{i=1}^N$.
\end{enumerate}
From  \eqref{setp26} it follows that 
the matrix $U\in \mathcal U_N$ has upper quasiseparable generators with orders
equal one.

All the given parameters define completely the matrix pair $(A,B)$ from the 
class $\mathcal P_N$. 

Each step of structured QZ should update these parameters 
while maintaining the minimal
orders of generators. However, structured algorithms for the multiplication 
of quasiseparable matrices may output redundant generators for the product 
matrix. For this reason, we will need an algorithm that compresses generators
to minimal order. The algorithm we use derives from previous work and is described in detail in the Appendix.

\section{The QZ step via generators}
\subsection{Classical QZ}
\setcounter{equation}{0}
Let $(A,B)$ be a pair of $N\times N$ matrices, with $A=(a_{ij})_{i,j=1}^N$
 upper Hessenberg and $B=(b_{ij})_{i,j=1}^N$ upper triangular. 
 The implicit QZ step  applied to $(A,B)$ consists in 
the computation of
unitary matrices $Q$ and $Z$ such that
\be\label{sten21}
\begin{array}{cc}
\mbox{the matrix}\;A_1=Q^*AZ\;\mbox{is upper Hessenberg}; \\
\mbox{the matrix}\;B_1=Q^*BZ\;\mbox{is upper triangular};
\end{array}
\end{equation}
and,  in addition,   some initial conditions are satisfied.
In the case of the single--shift implicit QZ step, the unitary matrices $Q$ and $Z$ are upper Hessenberg and take the form  
\be\label{smap22}
Q=\tl Q_1\tl Q_2\cdots\tl Q_{N-1},\quad Z=\tl Z_1\tl Z_2\cdots\tl Z_{N-1},
\end{equation}
where
\be\label{smip22}
\tl Q_i=I_{i-1}\oplus Q_i\oplus I_{N-i-1},\quad
\tl Z_i=I_{i-1}\oplus Z_i\oplus I_{N-i-1}
\end{equation}
and $Q_i,Z_i$ are complex Givens rotation matrices. The first Givens matrix $Q_1$ is chosen so that
\be\label{ssh1}
Q_1^*\left(\ba{c}a_{11}-\al\\a_{21}\ea\right)=
\left(\ba{c}\ast\\0\ea\right),
\end{equation}
where the shift $\al\in \mathbb C$ is an approximation of the eigenvalue that is currently being computed. The next Givens matrices $Q_i$ and $Z_i$ are computed and applied to the pencil using a bulge-chasing technique, which we show graphically on a $4\times 4$ example (see e.g., \cite{GVL}):

$$
\left(\begin{array}{cccc}
\times&\times&\times&\times\\
\times&\times&\times&\times\\
0&\times&\times&\times\\
0&0&\times&\times
\end{array}\right),
\left(\begin{array}{cccc}
\times&\times&\times&\times\\
0&\times&\times&\times\\
0&0&\times&\times\\
0&0&0&\times
\end{array}\right)
\overset{Q_1^*\cdot}{\longrightarrow}
\left(\begin{array}{cccc}
\times&\times&\times&\times\\
\times&\times&\times&\times\\
0&\times&\times&\times\\
0&0&\times&\times
\end{array}\right),
\left(\begin{array}{cccc}
\times&\times&\times&\times\\
\times&\times&\times&\times\\
0&0&\times&\times\\
0&0&0&\times
\end{array}\right)
$$

$$
\overset{\cdot Z_1}{\longrightarrow}
\left(\begin{array}{cccc}
\times&\times&\times&\times\\
\times&\times&\times&\times\\
\times&\times&\times&\times\\
0&0&\times&\times
\end{array}\right),
\left(\begin{array}{cccc}
\times&\times&\times&\times\\
0&\times&\times&\times\\
0&0&\times&\times\\
0&0&0&\times
\end{array}\right)
\overset{Q_2^*\cdot}{\longrightarrow}
\left(\begin{array}{cccc}
\times&\times&\times&\times\\
\times&\times&\times&\times\\
0&\times&\times&\times\\
0&0&\times&\times
\end{array}\right),
\left(\begin{array}{cccc}
\times&\times&\times&\times\\
0&\times&\times&\times\\
0&\times&\times&\times\\
0&0&0&\times
\end{array}\right)
$$

$$
\overset{\cdot Z_2}{\longrightarrow}
\left(\begin{array}{cccc}
\times&\times&\times&\times\\
\times&\times&\times&\times\\
0&\times&\times&\times\\
0&\times&\times&\times
\end{array}\right),
\left(\begin{array}{cccc}
\times&\times&\times&\times\\
0&\times&\times&\times\\
0&0&\times&\times\\
0&0&0&\times
\end{array}\right)
\overset{Q_3^*\cdot}{\longrightarrow}
\left(\begin{array}{cccc}
\times&\times&\times&\times\\
\times&\times&\times&\times\\
0&\times&\times&\times\\
0&0&\times&\times
\end{array}\right),
\left(\begin{array}{cccc}
\times&\times&\times&\times\\
0&\times&\times&\times\\
0&0&\times&\times\\
0&0&\times&\times
\end{array}\right)
$$

$$
\overset{\cdot Z_3}{\longrightarrow}
\left(\begin{array}{cccc}
\times&\times&\times&\times\\
\times&\times&\times&\times\\
0&\times&\times&\times\\
0&0&\times&\times
\end{array}\right),
\left(\begin{array}{cccc}
\times&\times&\times&\times\\
0&\times&\times&\times\\
0&0&\times&\times\\
0&0&0&\times
\end{array}\right)
$$

We provide for further reference a pseudocode description of a bulge-chasing sweep:
\begin{enumerate}
\item
determine $Q_1^*$ from \eqref{ssh1} and compute $A\leftarrow Q_1^*A$ 
\item
{\bf for } $k=1:N-2$ {\bf do}
\item
$\qquad$ $B\leftarrow Q_k^*B$
\item
$\qquad$ determine $Z_k$ that eliminates bulge in $B$
\item
$\qquad$ $A\leftarrow A*Z_k$
\item
$\qquad$ $B\leftarrow B*Z_k$
\item
$\qquad$ determine $Q_{k+1}^*$ that eliminates bulge in $A$
\item
$\qquad$ $A\leftarrow Q_{k+1}^*A$
\item 
{\bf end } {\bf do}
\item
$B\leftarrow Q_{N-1}^*B$
\item
determine $Z_{N-1}$ that eliminates bulge in $B$
\item
$A\leftarrow A*Z_{N-1}$
\item
$B\leftarrow B*Z_{N-1}$.
\end{enumerate}

\subsection{Structured QZ}
We present now a fast  adaptation  of the implicit 
single-shift QZ algorithm for  an input matrix  pair  $(A, B)\in\mathcal P_N$.
The modified algorithm works on the generators of the two  matrices  and this 
explains why it is referred to as a   structured  implicit QZ iteration with 
single shift. More specifically, the input pair $(A,B)=(A_0, B_0)$ is  first 
represented by means of a linear set of generators as explained  in 
Subsection \ref{psec}.  Then a high--level description of structured QZ 
iteration goes as follows.
 
\begin{enumerate}
\item Given the number $\alpha\in \mathbb C$, perform one step of the 
implicit single-shift QZ algorithm  \eqref{qzss} by computing subdiagonal 
entries of the matrix $A_1$, diagonal entries of the matrix $B_1$, vectors of
perturbation $\B z_1,\B w_1,\B p_1,\B q_1$ as well as upper triangular 
generators of the matrix $V_1$ and upper quasiseparable generators of the 
matrix $U_1$ (with redundant orders).
\item Compress the representations of $V_1$ and $U_1$ by using 
the compression algorithm in the Appendix.
\end{enumerate}
\medskip

Let us see how the QZ computation
 at step 1  above is efficiently performed by
working on the set of generators. 
When appropriate, we will reference lines from the classical QZ pseudocode shown above. 

\vspace{10pt}

{\bf Structured QZ bulge-chasing algorithm}

{\bf Input:} the shift $\alpha \in \mathbb C$,
subdiagonal entries $\sg^A_k$ for the matrix $A$, upper triangular
generators $g_V(k)$, $h_V(k)$, $b_V(k)$ for the matrix $V$, diagonal entries 
$d_B(k)$ for the matrix $B$, upper quasiseparable generators $g_U(k)$, 
$h_U( k)$, $b_U(k)$ for the matrix $U$, and perturbation vectors  ${\bf p}$, 
${\bf q}$, ${\bf z}$ and ${\bf w}$.

{\bf Output:} subdiagonal entries $\sg^{A_1}_k$ for the matrix $A_1$, upper 
triangular generators $g_{V_1}(k)$, $h_{V_1}(k)$, $b_{V_1}(k)$ of redundant orders for 
the matrix $V_1$, diagonal entries $d_{B_1}(k)$ for the matrix $B_1$, 
upper quasiseparable generators $g_{U_1}(k)$, $h_{U_1}(k)$, $b_{U_1}(k)$  of redundant orders
for the matrix $U_1$,  perturbation vectors  ${\bf p_1}$, ${\bf q_1}$, 
${\bf z_1}$, ${\bf w_1}$, and, if needed, the matrices $Q_k$ and $Z_k$.

\vspace{10pt}

\begin{enumerate}
\item Compute  
\be\nonumber
\sg_1^V=\sg_1^A+z(2)w^*(1),\quad \epsilon_1=g_V(1)h_V(1)-z(1)w^*(1).
\end{equation}
{\em (Pseudocode line 1).} Determine a complex Givens transformation matrix $Q_1$ from the condition
\be\nonumber
Q_1^*\left(\ba{c}\epsilon_1-\alpha\\\sg^A_1\ea\right)=
\left(\ba{c}\ast\\0\ea\right).
\end{equation}
Compute 
\be\nonumber
\Gamma_2=
Q_1^*\left(\ba{cc}g_V(1)h_V(1)&g_V(1)b_V(1)\\\sg^V_1&g_V(2)\ea\right).
\end{equation}
and determine the matrices $\tl g_V(2),\bt^V_2,f^V_2,\phi^V_2$ of sizes
$1\times(r^V_2+1),1\times(r^V_2+1),1\times1,1\times r^V_2$ from the partitions
\be\nonumber
\Gamma_2=\left[\ba{c}\tl g_V(2)\\\bt^V_2\ea\right],\quad
\beta^V_2=\left(\ba{cc}f^V_2&\phi^V_2\ea\right).
\end{equation}
Compute
\be\nonumber
\left(\ba{c}z^{(1)}(1)\\\chi_2\ea\right)=Q_1^*\left(\ba{c}z(1)\\z(2)\ea\right)
\end{equation}
with the numbers $z^{(1)}(1),\chi_2$.

Compute
\be\nonumber
f^A_2=f^V_2-\chi_2w^*(1),\quad \varphi^A_2=\phi_2^V.
\end{equation}

\item Set
\be\nonumber
\gamma_1=w(1),\;c_1=p(1),\;\theta_1=q(1).
\end{equation}
Compute
\be\nonumber
d_U(1)=d_B(1)+p(1)q^*(1)
\end{equation}
and set
\be\nonumber
f^U_1=d_U(1),\;\phi^U_1=g_U(1),
\end{equation}
\be\nonumber
f^B_1=d_B(1),\quad\varphi_1=g_U(1).
\end{equation}
\item  For $k=1,\dots,N-2$ perform the following.
\begin{enumerate}
\item {\em (Pseudocode line 3).} Compute the numbers
\be\nonumber
\begin{gathered}
\varepsilon_k=\varphi_kh_U(k+1)-c_kq^*(k+1),\\
\epsilon_{k+1}=\varphi^A_kh_V(k+1)-\chi_kw^*(k+1),\\
d_U(k+1)=d_B(k+1)+p(k+1)q^*(k+1)
\end{gathered}
\end{equation}
and the $2\times2$ matrix $\Phi_k$ by the formula
\be\nonumber
\Phi_k=Q_k^*\left(\ba{cc}f^B_k&\varepsilon_k\\0&d_B(k+1)\ea\right).
\end{equation}
\item  {\em (Pseudocode line 4).} Determine a complex Givens rotation matrix $Z_k$ such that
\be\nonumber
\Phi_k(2,:)Z_k=\left(\ba{cc}0&\ast\ea\right).
\end{equation}
 Compute the $2\times2$ matrix $\Om_k$ by the formula
\be\nonumber
\Om_k=\left(\ba{cc}f^A_{k+1}&\epsilon_{k+1}\\0&\sg^A_{k+1}\ea\right)Z_k
\end{equation}
\item   {\em (Pseudocode line 7).} Determine a complex Givens rotation matrix $Q_{k+1}$ and the number
 $\sg^{A_1}_k$ such that
\be\nonumber
 Q^*_{k+1}\Om_k(:,1)=\left(\ba{c}\sg^{A_1}_k\\0\ea\right).
 \end{equation}
\item  {\em (Pseudocode lines 5 and 8).} Compute
\be\nonumber
\begin{gathered}
\Gamma'_{k+2}=Q^*_{k+1}\left(\ba{ccc}f^V_{k+1}&\phi_{k+1}^Vh_V(k+1)&
\phi_{k+1}^Vb_V(k+1)\\z(k+2)\gamma^*_k&\sg^V_{k+1}&g_V(k+2)\ea\right),\\
\Gamma_{k+2}=\Gamma'_{k+2}\left(\ba{cc}Z_k&0\\0&I_{r^V_{k+2}}\ea\right),
\end{gathered}
\end{equation}
\be\nonumber
C_{k+2}=\left(\ba{ccc}1&0&0\\0&h_V(k+1)&b_V(k+1)\ea\right)
\left(\ba{cc}Z_k&0\\0&I_{r^V_{k+2}}\ea\right).
\end{equation}
and determine the matrices $\tl d_{V_1}(k+2),\tl g_{V_1}(k+2),\bt^V_{k+2},
f^V_{k+2},\phi^V_{k+2}$ of sizes 
$1\times1,1\times(r^V_{k+2}+1),1\times(r^V_{k+2}+1),
1\times1,1\times r^V_{k+2}$ from the partitions
\be\nonumber
\Gamma_{k+2}=\left[\ba{cc}\tl d_{V_1}(k+2)&\tl g_{V_1}(k+2)\\
\ast&\bt^V_{k+2}\ea\right],\quad
\bt^V_{k+2}=\left[\ba{cc}f_{k+2}^V&\phi_{k+2}^V\ea\right]
\end{equation}
and the matrices $\tl h_{V_1}(k+2),\tl b_{V_1}(k+2)$ of sizes 
$(r^V_{k+1}+1)\times1,(r^V_{k+1}+1)\times(r^V_{k+2}+1)$ from the partition
\be\nonumber
C_{k+2}=\left(\ba{cc}\tl h_{V_1}(k+2)&\tl b_{V_1}(k+2)\ea\right)
\end{equation}
\item  {\em (Pseudocode lines 3 and 6).} Compute
\be\nonumber
\begin{gathered}
\hskip-0.2cm\La_{k+1}=\\Q_k^*
\left(\ba{ccc}f^U_k&\phi^U_kh_U(k+1)&\phi^U_kb_U(k+1)\\
p(k+1)\theta^*_k&d_U(k+1)&g_U(k+1)\ea\right)
\left(\ba{cc}Z_k&0\\0&I_{r^U_{k+1}}\ea\right),
\end{gathered}
\end{equation}
\be\nonumber
D_{k+1}=\left(\ba{ccc}1&0&0\\0&h_U(k+1)&b_U(k+1)\ea\right)
\left(\ba{cc}Z_k&0\\0&I_{r^U_{k+1}}\ea\right),
\end{equation}
and determine the matrices $\tl d_{U_1}(k),\tl g_{U_1}(k),\bt^U_{k+1},
f^U_{k+1},\phi^U_{k+1}$ of sizes $1\times1,1\times(r^U_{k+1}+1),
1\times(r^U_{k+1}+1),1\times1,1\times r^U_{k+1}$ from the partitions
\be\nonumber
\La_{k+1}=\left(\ba{cc}\tl d_{U_1}(k)&\tl g_{U_1}(k)\\
\ast&\bt^U_{k+1}\ea\right),\quad
\bt^U_{k+1}=\left(\ba{cc}f^U_{k+1}&\phi^U_{k+1}\ea\right)
\end{equation}
and the matrices $\tl h_{U_1}(k+1),\tl b_{U_1}(k+1)$ of sizes 
$(r^U_k+1)\times1,(r^U_k+1)\times(r^U_{k+1}+1)$ from the partition
\be\nonumber
D_{k+1}=\left(\ba{cc}\tl h_{U_1}(k+1)&\tl b_{U_1}(k+1)\ea\right).
\end{equation}
\item {\em (Update perturbation vectors for $A$).} Compute
\be\nonumber
\begin{gathered}
\left(\ba{c}z^{(1)}(k+1)\\\chi_{k+2}\ea\right)=
Q^*_{k+1}\left(\ba{c}\chi_k\\z(k+1)\ea\right),\\
\left(\ba{c}w^{(1)}(k)\\\g_{k+1}\ea\right)=
Z^*_k\left(\ba{c}\g_k\\w(k+1)\ea\right)
\end{gathered}
\end{equation}
with the numbers $z^{(1)}(k+1),w^{(1)}(k),\chi_{k+2},\g_{k+1}$.
\item  {\em (Update perturbation vectors for $B$).} Compute
\be\nonumber
\left(\ba{c}q^{(1)}(k)\\\theta_{k+1}\ea\right)=
Z_k^*\left(\ba{c}\theta_k\\q(k+1)\ea\right),\quad
\left(\ba{c}p^{(1)}(k)\\c_{k+1}\ea\right)=
Q_k^*\left(\ba{c}c_k\\p(k+1)\ea\right)
\end{equation}
with the numbers $q^{(1)}(k),p^{(1)}(k),\theta_{k+1}, c_{k+1}$.
\item Compute
\be\nonumber
f^A_{k+2}=f^V_{k+2}-\chi_{k+2}\gamma^*_{k+1},\quad 
\varphi^A_{k+2}=\phi_{k+2}^V.
\end{equation}
Compute
\be\nonumber
f^B_{k+1}=f^U_{k+1}-c_{k+1}\theta^*_{k+1},\quad\varphi_{k+1}=\phi^U_{k+1}.
\end{equation}
\end{enumerate}
\item {\em (Pseudocode line 10).} Compute the numbers
\be\nonumber
\varepsilon_{N-1}=\varphi_{N-1}h_U(N)-c_{N-1}q^*(N),\quad
d_U(N)=d_B(N)+p(N)q^*(N)
\end{equation}
and the $2\times2$ matrix $\Phi_{N-1}$ by the formula
\be\nonumber
\Phi_{N-1}=Q_{N-1}^*\left(\ba{cc}f^B_{N-1}&\varepsilon_{N-1}\\
0&d_B(N)\ea\right).
\end{equation}
\item {\em (Pseudocode line 11).} Determine a complex Givens rotation matrix $Z_{N-1}$ such that
\be\nonumber
\Phi_{N-1}(2,:)Z_{N-1}=\left(\ba{cc}0&\ast\ea\right).
\end{equation}
\item {\em (Pseudocode line 12.} Compute
\be\nonumber
\Gamma_{N+1}=\left(\ba{cc}f^V_N&\phi_N^Vh_V(N)\ea\right)Z_{N-1},
\end{equation}
\be\nonumber
C_{N+1}=\left(\ba{cc}1&0\\0&h_V(N)\ea\right)Z_{N-1},\quad \tl h^{(1)}_V(N)=1.
\end{equation}
and determine the numbers $\tl d_{V_1}(N+1),\tl g_{V_1}(N+1)$ from the 
partition
\be\nonumber
\Gamma_{N+1}=\left[\ba{cc}\tl d_{V_1}(N+1)&\tl g_{V_1}(N+1)\ea\right].
\end{equation} 
and $r^V_N+1$-dimensional columns $\tl h_{V_1}(N+1),\tl b_{V_1}(N+1)$ from the 
partition
\be\nonumber
C_{N+1}=\left(\ba{cc}\tl h_{V_1}(N+1)&\tl b_{V_1}(N+1)\ea\right).
\end{equation}
\item {\em (Pseudocode line 13).} Compute
\be\nonumber
\La_N=Q_{N-1}^*
\left(\ba{cc}f^U_{N-1}&\phi^U_{N-1}h_U(N)\\p(N)\theta^*_{N-1}&d_U(N)\ea\right)
Z_{N-1},
\end{equation}
\be\nonumber
D_N=\left(\ba{cc}1&0\\0&h_U(N)\ea\right)Z_{N-1},\quad \tl h^{(1)}_U(N)=1
\end{equation}
and determine the numbers $\tl d_{U_1}(N-1),\tl g_{U_1}(N-1),
\tl d_{U_1}(N)$ from the partition
\be\nonumber
\La_N=\left[\ba{cc}\tl d_{U_1}(N-1)&\tl g_{U_1}(N-1)\\\ast&\tl d_{U_1}(N)
\ea\right]
\end{equation}
and $r^U_{N-1}+1$-dimensional columns $\tl h_{U_1}(N),\tl b_{U_1}(N)$ from the 
partition
\be\nonumber
D_N=\left(\ba{cc}\tl h_{U_1}(N)&\tl b_{U_1}(N)\ea\right).
\end{equation}
Set
\be\nonumber
\tl h_{V_1}(N+1)=1,\quad\tl h_{U_1}(N+1)=1.
\end{equation}
\item {\em (Update perturbation vectors).} Compute
\be\nonumber
\left(\ba{c}w^{(1)}(N-1)\\w^{(1)}(N)\ea\right)=
Q^*_{N-1}\left(\ba{c}\g_{N-1}\\w(N)\ea\right),
\end{equation}
\be\nonumber
\hskip-0.6cm\left(\ba{c}p^{(1)}(N-1)\\p^{(1)}(N)\ea\right)=
Q^*_{N-1}\left(\ba{c}c_{N-1}\\p(N)\ea\right),\;
\left(\ba{c}q^{(1)}(N-1)\\q^{(1)}(N)\ea\right)=
Z^*_{N-1}\left(\ba{c}\theta_{N-1}\\q(N)\ea\right).
\end{equation}
Set
\be\nonumber
z^{(1)}(N)=\chi_N
\end{equation}
and compute
\be\nonumber
\sg^{A_1}_{N-1}=\tl d_{V_1}(N+1)-z^{(1)}(N)(w^{(1)}(N))^*.
\end{equation}
\item Compute
\be\nonumber
d^{(1)}_B(k)=d_{U_1}(k+1)-p^{(1)}(k)(q^{(1)}(k))^*.
\end{equation}

\item
{\em (Adjust indices).} Set
\be\nonumber
\begin{gathered}
g_{V_1}(k)=\tl g_{V_1}(k+1),\;k=1,\dots,N,\quad
h_{V_1}(k)=\tl h_{V_1}(k+1),\;k=2,\dots,N+1,\\
b_{V_1}(k)=\tl b_{V_1}(k+1),\;k=2,\dots,N;\\
d_{V_1}(1),d_{V_1}(N+1)\;\mbox{to be the $1\times0$ and $0\times1$ empty 
matrices},\\ d_{V_1}(k)=\tl d_{V_1}(k+1),\;k=2,\dots,N
\end{gathered}
\end{equation}
and
\be\nonumber
\begin{gathered}
g_{U_1}(k)=\tl g_{U_1}(k+1),\;k=1,\dots,N,\\
h_{U_1}(k)=\tl h_{U_1}(k+1),\;k=2,\dots,N+1,\\
b_{U_1}(k)=\tl b_{U_1}(k+1),\;k=2,\dots,N;\\
d_{U_1}(k)=\tl d_{U_1}(k+1),\;k=1,\dots,N.
\end{gathered}
\end{equation}
\end{enumerate}

Although the formal proof of this algorithm is done via multiplication algorithms for matrices with quasiseparable representations (see \cite[Theorem 31.4 and Theorem 36.4]{EGH2}), the explanation can also be given via the bulge chasing process used in the standard QZ algorithms. 

For the most part, the structured update process is carried out separately on $V$, $U$, $z$ and $w$. However, note the following correspondences between some quantities computed in the structured algorithm above and some quantities used in the classical approach:
\begin{itemize}
\item 
The matrices $\Phi _k$ at step 3a correspond to the bulges created in the matrix $B$ and are explicitly computed in order to determine the Givens matrices $Z_k$;
\item 
Analogously, the matrices $\Omega _k$ at step 3b correspond to the bulges created in the matrix $A$ and are explicitly computed in order to determine the Givens matrices $Q_{k+1}$;
\item
The diagonal entries of $B$ and the subdiagonal entries of $A$ are explicitly stored in $d_B$ and $\sigma_A$, respectively. The diagonal entries of $A$, on the other hand, can be easily computed from the generators. For instance,  the quantity $\epsilon_1=A(1,1)$ is computed at  step 1 in order to determine the first Givens matrix $Q_1$ that encodes the shift. Observe that the explicit computation of diagonal and subdiagonal entries is also crucial when performing deflation.
\end{itemize} 

\subsection{Complexity}\label{sec_complexity} 
A complexity estimate on the structured QZ algorithm above, applied to an $N\times N$ pencil and followed by the compression algorithm found in the Appendix, yields:
\begin{itemize}
\item
$143+4\eta+(N-2)(251+7\eta)$ floating-point operations for the structured QZ update,
\item
$33+(N-2)(3\eta+144)$ operations for the compression of $V$,
\item
$3+(N-2)(2\eta+30)$ operations for the compression of $U$,
\end{itemize}
which gives a total count of $179+4\eta+(N-2)(425+12\eta)$ operations per iteration. Here $\eta$ denotes the computational cost required to compute and apply a $2\times 2$ Givens matrix.

\section{Backward error analysis for companion pencils}\label{sec_error}

This section is concerned with the  backward error analysis of the structured QZ algorithm for companion pencils.
The backward stability of structured QR/QZ algorithms based on generator representations has not received much attention
in the literature so far. To our knowledge the only algorithm which  is shown to be provably backward  stable is the modification of the QR
iteration for a rank-one correction of Hermitian matrices presented in \cite{EBG}.

This lack of works on stability analysis is probably due to the fact that generator computations are rather involved and do not admit a compact representation in
terms of matrix manipulations. To circumvent this difficulty here we adopt the following hybrid point of view. We assume that the fast structured QR algorithm is backward stable. In this section based on  a first order perturbation analysis we show  that the initial perturbed problem has the
same structure as the given problem and hence we derive the corresponding structured backward error expressed in terms of perturbations of the coefficients of the polynomial whose zeros are the computed roots.  The results of a thorough numerical/experimental investigation  are reported in the next section. These results confirm that the computed {\em a posteriori} bounds on the coefficients fit our analysis
whenever we suppose that the backward error introduced by the fast variant of the QZ iteration is a small multiple of the
same error for the customary algorithm. This is also in accordance with the conclusion stated in \cite{EBG}. 


The problem of determining whether the QR method applied to the Frobenius companion matrix yields a backward stable rootfinder is 
 examined in \cite{EM}, with a positive answer (Theorem 2.1) if the norm of the polynomial is small. Given a monic polynomial $p(x)=\sum_{j=0}^{n}a_j x^j$, let $A$ be its companion matrix and $E$ the perturbation matrix that measures the backward error of the QR method applied to $A$. Edelman and Murakami show that the coefficients of $s^{k-1}$ in the polynomial ${\rm det}(sI_n-A-E)-{\rm det}(sI_n-A)$ are given at first order by the expression
\begin{equation}\label{EM_formula}
\sum_{m=0}^{k-1}a_m \sum_{i=k+1}^n E_{i,i+m-k}-\sum_{m=k}^{n}a_m \sum_{i=1}^k E_{i,i+m-k},
\end{equation}
with $a_n=1$. 

A similar property holds for QZ applied to companion pencils. This result was proven in \cite{VD} by using a block elimination process, and later in \cite{LV03} and \cite{EK} via a geometric approach. See also \cite{LV06} for a backward error analysis that takes into account the effects of balancing.

In this section we rely on the results in \cite{EM} to derive explicit formulas for the polynomial backward error, at least for the case of nonsingular $B$. Here $a_n$ is no longer necessarily equal to $1$, although it is different from $0$. 

We consider companion pencils of the form
\begin{equation}\label{eq_comp_pencil}
A=\left(\begin{array}{cccc}
&&&-a_0\\
1&&&-a_1\\
&\ddots&&\vdots\\
&&1&-a_{n-1}
\end{array}\right),\quad
B=\left(\begin{array}{cccc}
1&&&\\
&\ddots&&\\
&&1&\\
&&&a_n
\end{array}\right)
\end{equation}
and we will denote the perturbation pencil as $(E,G)$, where $E$ and $G$ are matrices of small norm.

First, observe that \eqref{EM_formula} can be generalized to the computation of ${\rm det} (s(I_n+G)-(A+E))$. Indeed, at first order we have
\begin{eqnarray*}
&& s(I_n+G)-(A+E)=(I_n+G)[s I_n -(I_n+G)^{-1}(A+E)]\\
&& = (I_n+G)[s I_n -(I_n-G)(A+E)]\\
&& = (I_n+G)[sI_n-A-E+GA]
\end{eqnarray*}
and therefore
\begin{equation}\label{EM_generalized}
{\rm det} (s(I_n+G)-(A+E))={\rm det}(I_n+G){\rm det} (sI_n-A-E+GA),
\end{equation}
where (again at first order)
\begin{equation}
{\rm det}(I_n+G)= 1+{\rm tr}(G)
\end{equation}
and ${\rm det} (sI_n-A-E+GA)$ can be computed by applying \eqref{EM_formula} to
the pencil $sI_n-A$ with the perturbation matrix $E-GA$. In order to fix the 
notation, let us write at first order
\begin{equation}\label{det_approx}
\begin{gathered}
 {\rm det} (sI_n-A-E+GA)=
{\rm det} (sI_n-A)+\left(\sum_{j=0}^{n-1} (\Delta a)_j s^j\right).
\end{gathered}
 \end{equation}
 
 Now, recall that $B$ can be seen as a rank-one perturbation of the identity matrix: we can write $B=I_n+(a_n-1)e_ne_n^T$, with the usual notation $e_n=[0,\ldots ,0,1]^T$. So, the perturbed companion pencil is
\begin{equation}\label{pert_pencil}
s(B+G)-(A+E)=s(I_n+G)-(A+E)+s(a_n-1)e_ne_n^T.
\end{equation}
The Sherman-Morrison determinant formula applied to (\ref{pert_pencil}) gives
\begin{equation}\label{pert_pencil2}
\begin{gathered}
{\rm det} (s(B+G)-(A+E))=\\
{\rm det} (s(I_n+G)-(A+E))(1+s(a_n-1)e_n^T[s(I_n+G)-(A+E)]^{-1}e_n).
\end{gathered}
\end{equation}
We know how to compute the determinant in the right-hand side of 
(\ref{pert_pencil2}) thanks to the discussion above. Now we want to compute the
second factor. What we need is the $(n,n)$ entry of the matrix 
$[s(I_n+G)-(A+E)]^{-1}$, which can be written  as
\begin{eqnarray}
&&[s(I_n+G)-(A+E)]^{-1}(n,n) = \nonumber\\
&&={\rm det} (s(I_n+G)-(A+E))^{-1} \cdot{\rm det}(s(I_{n-1}+\tilde{G})-(\tilde{A}+\tilde{E})),\label{entry11}
\end{eqnarray}
where $\tilde{A}=A(1:n-1,1:n-1)$, $\tilde{E}=E(1:n-1,1:n-1)$ and $\tilde{G}=G(1:n-1,1:n-1)$. Observe that
 $\tilde{A}$ is a companion matrix for the polynomial $s^{n-1}$, so, using \eqref{EM_generalized}, we can write 
$${\rm det}(s(I_{n-1}+\tilde{G})-(\tilde{A}+\tilde{E}))={\rm det}(I_{n-1}+\tilde{G})\left[s^{n-1}+\sum_{j=0}^{n-2} (\Delta \tilde{a})_j s^j\right],$$
 where the coefficients $(\Delta \tilde{a})_j$
 can be computed from \eqref{EM_formula} with the perturbation matrix $\tilde{E}-\tilde{G}\tilde{A}$.
 
 From \eqref{EM_generalized}, \eqref{det_approx}, \eqref{pert_pencil2} and \eqref{entry11} we obtain:
\begin{eqnarray*}
&&{\rm det} (s(B+G)-(A+E))=\\
&&={\rm det}(s(I_n+G)-(A+E))+s(a_n-1){\rm det}(s(I_{n-1}+\tilde{G})-(\tilde{A}+\tilde{E}))=\\
&&={\rm det}(sI_n-A)+\sum_{j=0}^{n-1} (\Delta a)_j s^j+{\rm tr}(G){\rm det}(sI_n-A)+\\
&&+s(a_n-1)\left[\sum_{j=0}^{n-2} (\Delta \tilde{a})_j s^j+(1+{\rm tr}(\tilde{G}))s^{n-1}\right]=\\
&&=p(s)+{\rm tr}(G)q(s)+\sum_{j=0}^{n-1} (\Delta a)_j s^j+s(a_n-1)\left[\sum_{j=0}^{n-2} (\Delta \tilde{a})_j s^j+{\rm tr}(\tilde G)s^{n-1}\right],
\end{eqnarray*}
where $q(x)=x^n+\sum_{j=0}^{n-1}a_jx^j$. 

Finally, from ${\rm det} (s(B+G)-(A+E))$ we subtract the quantity ${\rm det}(sB-A)=p(s)$, and we obtain the following result:
\begin{prop}
With the above notation, the following equality is correct at first order:
\begin{eqnarray}
&&{\rm det} (s(B+G)-(A+E))-{\rm det}(sB-A)=\nonumber\\
&&={\rm tr}(G)q(s)+\sum_{j=0}^{n-1} (\Delta a)_j s^j+s(a_n-1)\left[\sum_{j=0}^{n-2} (\Delta \tilde{a})_j s^j+{\rm tr}(\tilde G)s^{n-1}\right].\label{BEG_generalization}
\end{eqnarray}
\end{prop}

\begin{remark}
A first-order estimate of ${\rm det} (A+E-s(B+G))-{\rm det}(A-sB)$ can also be obtained via a geometric approach, as shown in \cite{EM} for companion matrices and in \cite{LV03}, \cite{JV} and \cite{EK} for pencils. Note that this approach is related to the Leverrier rootfinding method: see \cite{mertzios} for an adaptation to the matrix pencil case.
\end{remark}

\section{Numerical results}

In this section we test the performance of the fast QZ method presented above and compare it to classical unstructured QZ applied to the companion pencil and to classical QR applied to the companion matrix. We have implemented our method in Matlab and in Fortran 90/95; both implementations are available online.\footnote{http://www.unilim.fr/pages\_perso/paola.boito/software.html}

The first set of experiments are performed in Matlab: we use the commands {\tt roots} for classical QR and {\tt eig} for classical QZ, whereas fast QZ is applied using our Matlab implementation of the algorithm described above. 

 In particular, we report absolute forward and backward errors, measured in $\infty$-norm. For each polynomial $p(x)$ of degree $N$, the forward error  is computed as 
 $$
 {\rm forward}\,{\rm error}=\max_{k=1,\dots,N}\min_{j=1,\dots,N}|\lambda_j-\alpha_k|,
 $$
 where the $\lambda_j$'s are the roots of $p(x)$ as determined by the eigensolver that is being studied, and the $\alpha_k$'s are the "exact" roots of $p(x)$, either known in advance or computed in high precision using Matlab's Symbolic Toolbox, unless otherwise specified. The backward error is computed as
 $$
 {\rm backward}\,{\rm error}=\max_{k=0,\dots,N}|\tilde{p}_k-p_k|,
 $$
where the $p_k$'s are the exact coefficients of $p(x)$, either known in advance or computed in high precision from the known exact roots, and the $\tilde{p}_k$'s are the coefficients computed in high precision from the $\lambda_j$'s. 
 
Our aim here is to provide experimental evidence pointing to the stability of our structured QZ method and to the better accuracy of QZ versus QR on some classes of polynomials. We do not report timings for the various methods, since the running times for our Matlab implementation cannot be compared to the running times of a built-in function such as {\tt eig} or {\tt roots}. 

Observe that the normalization of the polynomials is a crucial step for the proper functioning of QZ. Unlike the companion matrix, which is necessarily associated with a monic polynomial, the companion pencil allows for an arbitrary scaling of the polynomial. Unless otherwise specified, we normalize w.r.t.~the 2-norm of the vector of coefficients. The polynomials obtained from computed roots, used for computation of backward errors, are also normalized in 2-norm. For the purpose of computing backward errors, we have also experimented with normalization in a least-squares sense, as suggested in \cite{JV}, but in our examples we have generally found little difference between the 2-norm and the least-squares approach.

\vspace{10pt}

\begin{example}\label{ex_random}
Random polynomials.
\end{example}
We apply our structured QZ method to polynomials whose coefficients are random complex numbers with real and imaginary parts uniformly chosen in $[-1,1]$. Here $N$ denotes the degree. Such test polynomials are generally well conditioned and they are useful to study the behavior of our method as the polynomial degree grows larger. 

Table \ref{test_random} shows absolute forward  errors w.r.t.~the roots computed by the Matlab command {\tt roots}, as well as the average number of iterations per eigenvalue, which is consistent with the expected number of operations for the classical QZ method. Backward errors (not shown in the table) are of the order of the machine epsilon. For each degree, errors and the number of iterations are averaged over $10$ random polynomials.

\begin{table}[h!]
\begin{center}
\begin{tabular}{|c|c|c|}
\hline
$N$  & abs. forward error  & average n. iterations\\
\hline
$50$&$2.60$e$-14$&3.71\\
$100$&$1.30$e$-13$&$3.59$\\
$150$&$4.52$e$-13$&$3.45$\\
$200$&$5.90$e$-13$&$3.38$\\
$300$&$1.69$e$-12$&$3.28$\\
$400$&$4.34$e$-12$&$3.22$\\
$500$&$6.11$e$-12$&$3.18$\\
\hline
\end{tabular}
\end{center}
\caption{Errors and number of iterations for structured QZ applied to random polynomials: see Example \ref{ex_random}.}\label{test_random}
\end{table}

\begin{example}\label{ex_cyclotomic}
Cyclotomic polynomials $z^N-i$.
\end{example} 
This is another set of well-conditioned polynomials for which we consider degrees up to $500$. Table \ref{test_cyclotomic} shows forward errors with respect to the roots computed by the Matlab command {\tt roots}, as well as the average number of iterations per eigenvalue. Backward errors (not shown in the table) are of the order of the machine epsilon.

\begin{table}[h!]
\begin{center}
\begin{tabular}{|c|c|c|}
\hline
$N$  & abs.  forward error & average n. iterations\\
\hline
$50$&$1.01$e$-14$&$4.16$\\
$100$&$1.46$e$-14$&$4.01$\\
$150$&$3.55$e$-14$&$3.85$\\
$200$&$3.85$e$-14$&$3.88$\\
$300$&$7.23$e$-14$&$3.69$\\
$400$&$1.66$e$-13$&$3.66$\\
$500$&$1.78$e$-13$&$3.75$\\
\hline
\end{tabular}
\end{center}
\caption{Errors and number of iterations for structured QZ applied to cyclotomic polynomials: see Example \ref{ex_cyclotomic}.}\label{test_cyclotomic}
\end{table}

\begin{example}\label{ex_deg20}
\end{example}
We consider here a few examples taken from \cite{EM}:
\begin{enumerate}
\item $p(x)=x^{20}+x^{19}+\dots +x+1$,
\item the polynomial with roots equally spaced in the interval  $[-2.1,1.9]$,
\item the Chebyshev polynomial of first kind of degree $20$,
\item the Bernoulli polynomial of degree $20$.
\end{enumerate}
Table \ref{test_deg20} shows forward errors for fast and classical QZ (after normalization of the polynomial) and for classical balanced QR (before normalization), as well as the maximum eigenvalue condition number.
\begin{table}[h]
\begin{center}
\begin{tabular}{|c|c|c|c|c|}
\hline
polynomial  & fast QZ & classical QZ & QR & max condeig\\
\hline
$(1)$&$3.58$e$-15$&$1.09$e$-15$&$2.14$e$-15$&$1.38$\\
$(2)$&$7.87$e$-13$&$1.98$e$-13$&$5.24$e$-13$&$2.71$e$+4$\\
$(3)$&$4.16$e$-10$&$2.18$e$-11$&$3.62$e$-12$&$1.86$e$+296$\\
$(4)$&$4.00$e$-3$&$4.00$e$-3$&$4.00$e$-3$&$9.14$e$+6$\\
\hline
\end{tabular}
\end{center}
\caption{Forward errors for polynomials in Example \ref{ex_deg20}}\label{test_deg20}
\end{table}

Following \cite{EM}, we also test our backward error analysis, given in Section \ref{sec_error}, on these examples. Table \ref{be_deg20} shows the logarithm in base 10, rounded to an entire number,  of the computed  and of the predicted backward errors for each coefficient of each polynomial. The predicted error is computed by applying \eqref{BEG_generalization}. 

As for the choice of $E$ and $G$, we apply the backward error analysis given in \cite{EBG}, Theorem 4.1. This analysis implies that the backward error matrix for the QR method applied to a small rank perturbation of a Hermitian $N\times N$ matrix $M$ is, roughly speaking, bounded by a small multiple of $\varepsilon N^2 \|M\|_F$, where $\varepsilon$ is the machine epsilon. In view of this result, each nonzero entry of $E$ and $G$ can be taken as a small multiple of $N\varepsilon$, where $N$ is the degree of the polynomial. For the sparsity pattern of $E$ and $G$, we follow the ideas in \cite{EM}. So, $E$ and $G$ are chosen here as
  {\tt E=10*N*eps*triu(ones(N),-2)} and {\tt G=10*N*eps*triu(ones(N),-1)}, with {\tt N}=20, in order to model the backward error introduced by the QZ iterations. 
  Note that, under these assumptions, the experimental results confirm the theoretical analysis.

\begin{table}[h]
\begin{center}
\begin{tabular}{|l|c|c|c|c|c|}
\hline
&$(1)$  & $(2)$ & $(3)$ & $(4)$ & Ex. \ref{ex_QZvsQR} \\
\hline
$z^{20}$&$-14,-15,-13$&$-17,-17,-13$&$-16,-16,-13$&$-19,-19,-13$&$-26,-27,-13$\\
$z^{19}$&$-14,-15,-13$&$-17,-17,-15$&$-16,-16,-14$&$-18,-18,-16$&$-14,-15,-13$\\
$z^{18}$&$-14,-15,-13$&$-16,-16,-14$&$-15,-16,-13$&$-17,-18,-15$&$-15,-15,-14$\\
$z^{17}$&$-14,-15,-13$&$-16,-15,-14$&$-15,-16,-14$&$-17,-17,-16$&$-15,-15,-13$\\
$z^{16}$&$-15,-15,-13$&$-15,-16,-13$&$-15,-15,-12$&$-17,-17,-15$&$-15,-15,-13$\\
$z^{15}$&$-15,-15,-13$&$-15,-15,-14$&$-15,-15,-14$&$-16,-17,-15$&$-15,-15,-13$\\
$z^{14}$&$-14, -16,-13$&$-15,-15,-13$&$-15,-15,-12$&$-16,-16,-14$&$-15,-15,-13$\\
$z^{13}$&$-14,-15,-13$&$-15,-14,-14$&$-15,-16,-14$&$-15,-16,-15$&$-15,-15,-13$\\
$z^{12}$&$-15,-15,-13$&$-16,-15,-12$&$-15,-16,-12$&$-16,-16,-13$&$-15,-15,-13$\\
$z^{11}$&$-15,-15,-13$&$-15,-15,-14$&$-15,-15,-15$&$-15,-15,-14$&$-15,-15,-13$\\
$z^{10}$&$-15,-15,-13$&$-16,-16,-12$&$-15,-16,-13$&$-16,-16,-13$&$-15,-15,-13$\\
$z^{9}$&$-15,-15,-13$&$-15,-15,-14$&$-16,-16,-16$&$-14,-15,-14$&$-15,-15,-13$\\
$z^{8}$&$-15,-37,-13$&$-16,-15,-13$&$-16,-16,-14$&$-15,-16,-12$&$-15,-15,-13$\\
$z^{7}$&$-15,-15,-13$&$-15,-16,-15$&$-15,-16,-17$&$-14,-15,-14$&$-15,-15,-13$\\
$z^{6}$&$-14,-15,-13$&$-16,-15,-13$&$-15,-15,-15$&$-16,-40,-12$&$-15,-16,-13$\\
$z^{5}$&$-14,-16,-14$&$-16,-16,-16$&$-15,-16,-19$&$-14,-15,-14$&$-15,-15,-12$\\
$z^{4}$&$-14,-16,-14$&$-16,-16,-14$&$-15,-15,-16$&$-16,-16,-12$&$-15,-15,-13$\\
$z^{3}$&$-15,-15,-14$&$-17,-18,-17$&$-16,-16,-19$&$-15,-16,-17$&$-15,-14,-12$\\
$z^{2}$&$-15,-15,-28$&$-15,-17,-15$&$-15,-18,-18$&$-16,-17,-13$&$-15,-15,-13$\\
$z^{1}$&$-15,-15,-14$&$-18,-20,-17$&$-15,-$Inf,$-19$&$-16,-16,-15$&$-15,-15,-12$\\
$z^{0}$&$-15,-15,-13$&$-17,-19,-13$&$-47,-$Inf,$-13$&$-16,-17,-13$&$-15,-27,-13$\\
\hline
\end{tabular}
\end{center}
\caption{Backward errors for polynomials in Examples \ref{ex_deg20} and \ref{ex_QZvsQR}. Each entry of the table contains the logarithm in base 10 of the backward errors, in the following order: computed
via fast QZ, computed via classical QZ, predicted according to \eqref{BEG_generalization}.}\label{be_deg20}
\end{table}

\begin{example}\label{ex_QZvsQR}
QZ vs. QR. 
\end{example}
In this example with heavily unbalanced coefficients, the (classical or structured) QZ method applied to the companion pencil computes the roots with  better accuracy than QR applied to the companion matrix. We take here the polynomial $p=\sum_{k=0}^{20} p_kx^k$, where  $p_k=10^{6(-1)^{(k+1)} -3}$ for $k=0,\ldots,20$. Table \ref{test_QZvsQR} shows absolute forward and backward errors for several eigenvalue methods, with or without normalizing the polynomial before the computation of its roots. See Table \ref{be_deg20} for backward errors.

\begin{table}[h!]
\begin{center}
\begin{tabular}{|r|l|l|}
\hline
Method&Backward error&Forward error\\
\hline
Fast QZ & $2.45$e$-15$&$4.28$e$-15$\\
Classical QZ & $1.49$e$-15$&$3.22$e$-15$\\
Balanced QR after normalization& $1.22$e$-4$&$6.27$e$-9$\\
Unbalanced QR after normalization & $1.22$e$-4$ &$1.72$e$-15$\\
Balanced QR, no normalization & $1.22$e$-4$&$5.86$e$-9$\\
Unbalanced QR, no normalization & $1.22$e$-4$&$2.72$e$-15$\\
\hline
\end{tabular}
\end{center}
\caption{Errors for several versions of QZ and QR applied to a polynomial with unbalanced coefficients (Example \ref{ex_QZvsQR}).}
\label{test_QZvsQR}
\end{table}

\vspace{10pt}

\begin{table}[h]
\begin{center}
\begin{tabular}{|c|c|c|c|}
\hline
$N$  & abs. forward error & fast QZ time &LAPACK time \\
\hline
$50$&$1.73$e$-14$&$1.08$e$-2$&$7.10$e$-3$\\
$60$&$2.96$e$-14$&$1.66$e$-2$&$1.14$e$-2$\\
$70$&$3.32$e$-14$&$2.01$e$-2$&$1.76$e$-2$\\
$80$&$6.04$e$-14$&$2.43$e$-2$&$2.49$e$-2$\\
$90$&$7.74$e$-14$&$2.88$e$-2$&$3.51$e$-2$\\
$100$&$1.10$e$-13$&$3.26$e$-2$&$4.35$e$-2$\\
$150$&$1.60$e$-13$&$7.77$e$-2$&$1.35$e$-1$\\
$200$&$2.36$e$-13$&$1.29$e$-1$&$3.11$e$-1$\\
$300$&$8.44$e$-13$&$2.70$e$-1$&$9.91$e$-1$\\
$400$&$1.65$e$-12$&$4.66$e$-1$&$2.32$\\
$500$&$1.85$e$-12$&$7.27$e$-1$&$4.70$\\
\hline
\end{tabular}
\end{center}
\caption{Errors and running times (measured in seconds) for fast QZ implemented in Fortran versus the LAPACK implementation.}\label{table_timings}
\end{table}

\begin{figure}
\includegraphics[width=0.8\textwidth]{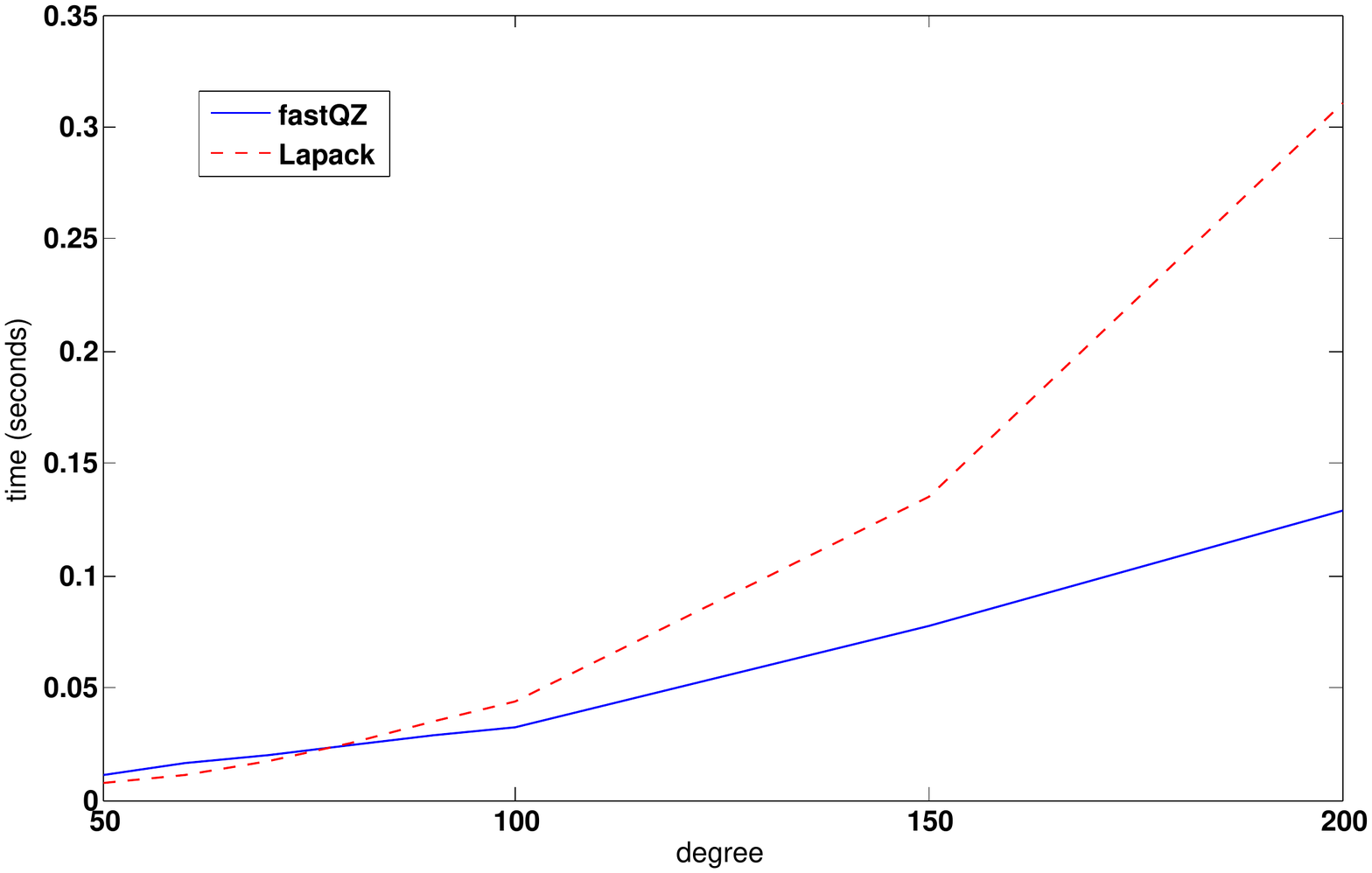}
\caption{Timings for our fast QZ algorithm (blue solid line) and for LAPACK (red dashed line).}\label{fig_comparison}
\end{figure}
\begin{figure}
\includegraphics[width=0.8\textwidth]{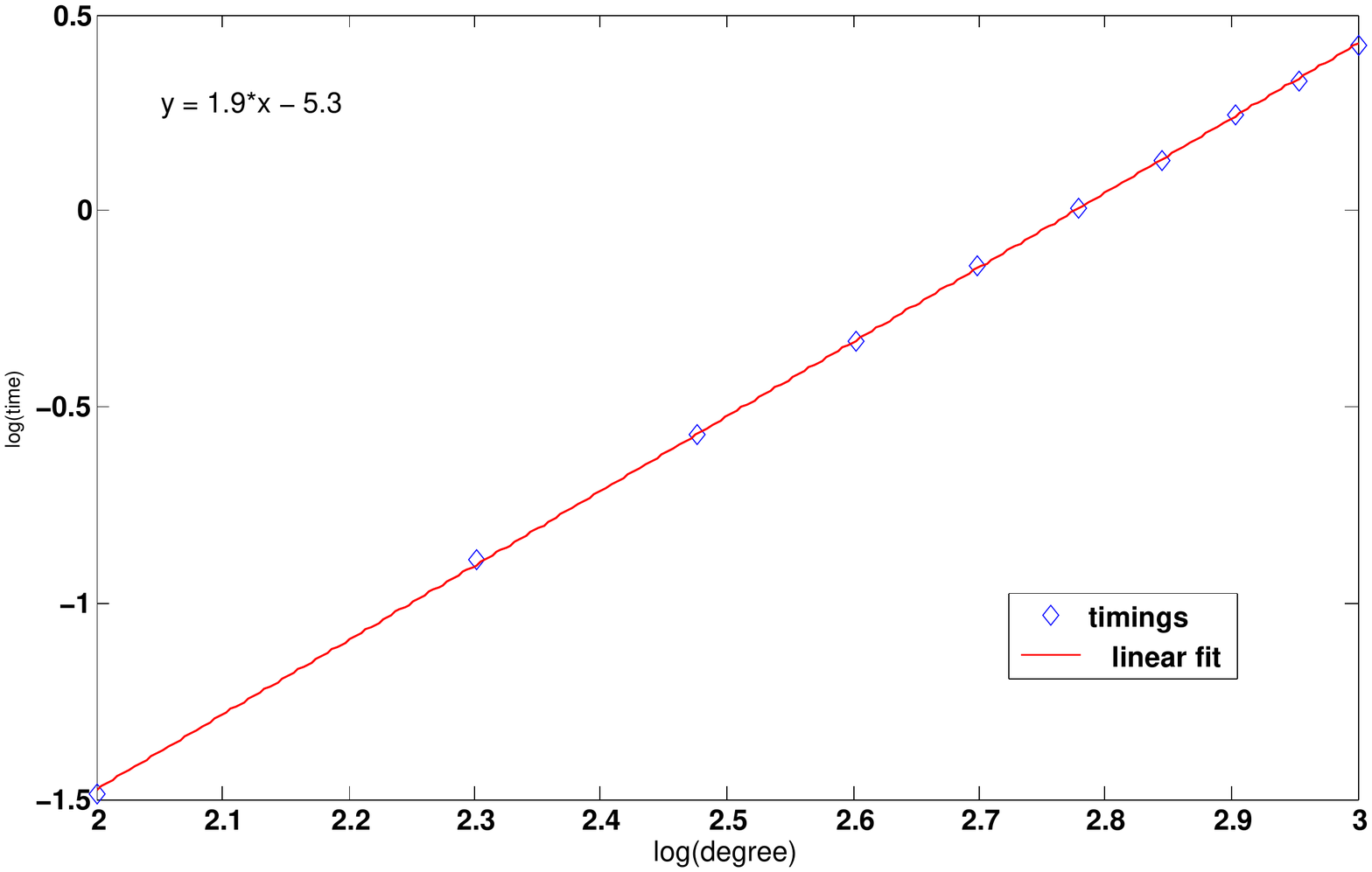}
\caption{Linear fit on log-log plot of running times. Polynomial degrees go from $100$ to $1000$.}\label{fig_asymptotic}
\end{figure}

We next show the results of experiments using the Fortran implementation of our method. The main goal of these experiments is to show that the structured method is actually faster than the classical LAPACK implementation for sufficiently high degrees, and that running times grow quadratically with the polynomial degree, as predicted by the complexity analysis. These experiments were performed on a MacBook Pro, using the {\tt gfortran} compiler. 

\begin{example}
Timings for random polynomials.
\end{example}
We take polynomials with random coefficients as in Example \ref{ex_random}. Table \ref{table_timings} shows forward absolute errors with respect to the roots computed by LAPACK (subroutine ZGEGV), as well as timings (measured in seconds). Results for each degree are averages over 10 trials. 

Figure \ref{fig_comparison} shows running times versus polynomial degrees for Lapack and for our structured implementation. The crossover point for this experiment appears to be located at degree about 80.

Figure \ref{fig_asymptotic} shows a log-log plot of running times for our implementation, together with a linear fit. The slope of the fit is $1.9$, which is consistent with the claim that the structured method has complexity $\mathcal{O}(n^2)$.

\section{Conclusions}

In this work we have presented and tested a new structured version of the single-shift, implicit QZ method. Our algorithm is designed for the fast computation of eigenvalues of matrix pencils belonging to a class $\mathcal P_N$ that includes companion and Lagrange pencils. The quasiseparable structure of such pencils allows us to achieve quadratic complexity. 

Our numerical experience says that the fast structured implementation of the QZ algorithm applied to a companion pencil provides a fast and provably backward stable root-finding method.

\appendix

\section{The compression of generators}\label{sec_compression}

\setcounter{equation}{0}

We present an algorithm that takes as input the possibly 
redundant quasiseparable generators of a unitary matrix and outputs generators 
of minimal order. This algorithm is justified in a similar way as the compression algorithm 
\cite[Theorem 7.5]{EGH1} but with computations in the forward direction.

Let $U=\{U_{ij}\}_{i,j=1}^N$ be a block unitary matrix with entries of
sizes $m_i\times n_j$, lower quasiseparable generators
$p(i)\;(i=2,\dots,N),\;q(j)\;(j=1,\dots,N-1),\;a(k)\;(k=2,\dots,N-1)$
of orders $r^L_k\;(k=1,\dots,N-1)$, upper quasiseparable generators
$g(i)\; (i=1,\dots,N-1),\;
h(j)\;(j=2,\dots,N),\;b(k)\;(k=2,\dots,N-1)$ of orders
$r^U_k\;(k=1,\dots,N-1)$ and diagonal entries
$d(k)\;(k=1,\dots,N)$. Set
\be\label{nnov8}
\begin{gathered}
\rho_0=0,\;\rho_k=\min\{n_k+\rho_{k-1},\;r^L_k\},\;k=1,\dots,N-1,\\
\nu_k=n_k+\rho_{k-1}-\rho_k,\;k=1,\dots,N-1,\quad\nu_N=n_N+\rho_{N-1},\\
s_0=0,\quad s_k=m_k+s_{k-1}-\nu_k,\quad k=1,\dots,N-1.
\end{gathered}
\end{equation}

Then all the numbers $s_k$ are nonnegative and the matrix $U$ has upper
quasiseparable generators of orders $s_k\;(k=1,\dots,N-1)$. A set of such upper
quasiseparable generators are obtained using the following algorithm. 

\vspace{10pt}

{\bf Compression algorithm}

{\bf Input:} lower quasiseparable generators $p(j)$, $q(j)$, $a(j)$ and upper 
quasiseparable generators $g(j)$, $h( j)$, $b(j)$ of possibly redundant orders 
for the matrix $U$.

{\bf Output:} upper quasiseparable generators $g_s(j)$, $h_s(j)$, $b_s(j)$ of 
minimal order for $U$.

\begin{enumerate}

\item
Set $X_0,Y_0,z_0$ to be the $0\times0$ empty matrices and
$p(1),a(1),h(1),b(1),h_s(1)$ to be empty matrices of sizes $m_1\times 0,
r^L_1\times 0,0\times n_1,0\times r^U_1,0\times n_1$ respectively.

\item
For $k=1,\dots,N-1$ perform the following. Determine an 
$(n_k+\rho_{k-1})\times(n_k+\rho_{k-1})$ unitary matrix $W_k$ and an
$r^L_k\times\rho_k$ matrix $X_k$ such that
\begin{equation}\nonumber
\left(\ba{cc}a(k)X_{k-1}&q(k)\end{array}\right)W_k^*=
\left(\ba{cc}0_{r^L_k\times\nu_k}&X_k\end{array}\right).
\end{equation}
$W_k$ can be computed, for instance, via the usual Givens or Householder methods.
\item
Compute the $(m_k+s_{k-1})\times(n_k+\rho_{k-1})$ matrix
\be\nonumber
Z_k=\left(\ba{cc}z_{k-1}&h_s(k)\\p(k)X_{k-1}q(k)&d(k)\ea\right)W_k^*.
\end{equation}
\item
Determine the matrices $\Theta_k,\Delta_k,h'_k,h''_k$ of sizes
$s_{k-1}\times\nu_k,m_k\times\nu_k,s_{k-1}\times\rho_k,m_k\times\rho_k$
from the partition
\be\nonumber
Z_k=\left[\ba{cc}\Theta_k&h'_k\\\Delta_k&h''_k\ea\right].
\end{equation}
Observe that the submatrix $\left(\ba{c}\Theta_k\\\Delta_k\ea\right)$ has 
orthonormal columns. 
\item
Determine an $(s_{k-1}+m_k)\times(s_{k-1}+m_k)$ unitary 
matrix $F_k$ from the condition
\be\nonumber
F^*_k\left(\ba{c}\Theta_k\\\Delta_k\ea\right)=
\left(\ba{c}I_{\nu_k}\\0_{s_k\times\nu_k}\ea\right).
\end{equation}
\item
Determine the matrices $h_F(k),d_F(k),b_s(k),g_s(k)$ of sizes
$s_{k-1}\times\nu_k,m_k\times\nu_k,s_{k-1}\times s_k,m_k\times s_k$ from the
partition
\be\nonumber
F_k=\left[\ba{cc}h_F(k)&b_s(k)\\d_F(k)&g_s(k)\ea\right].
\end{equation}
\item
Compute the matrices $Y_k$ of size $s_k\times r^U_k$ and $z_k$
of the size $s_k\times\rho_k$ by the formulas
\be\nonumber
Y_k=g^*_s(k)g(k)+b^*_s(k)Y_{k-1}b(k),\quad z_k=g^*_s(k)h''_k+b^*_s(k)h'_k.
\end{equation}
\item
Compute
\be\nonumber
h_s(k+1)=Y_kh(k+1).
\end{equation}
\item
Set
\be\nonumber
F_N=\left[\ba{cc}z_{N-1}&Y_{N-1}h(N)\\p(N)X_{N-1}&d(N)\ea\right].
\end{equation}
\end{enumerate}

For a matrix from the class ${\mathcal U}_N$ we have $m_i=n_i=1,\;i=1,\dots,N$
and $r^L_k=1,\;k=1,\dots,N-1$ which by virtue of (\ref{nnov8}) implies
$\rho_k=s_k=1,\;k=1,\dots,N-1$.

For a matrix from the class ${\mathcal V}_N$ we determine the sizes of blocks 
via (\ref{aprmn18s}) and the orders of lower generators
$r^L_k=1,\;k=1,\dots,N$. Hence using (\ref{nnov8}) we obtain
$s_1=1,\;s_k=2,\;k=2,\dots,N-1$.

\vspace{10pt}

\begin{remark}
With the above hypotheses, the matrix $U$ admits the (non-unique) factorization
\begin{equation}\nonumber
U=W\cdot F,
\end{equation}
where $W$ is a block lower triangular unitary matrix with block
entries of sizes $m_i\times\nu_j\;(i,j=1,\dots,N)$ and $F$ is a
block upper triangular unitary matrix with block entries of sizes
$\nu_i\times n_j\;(i,j=1,\dots,N)$. Moreover one can choose the
matrix $W$ in the form
\begin{displaymath}
W=\tl W_{N-1}\cdots\tl W_1,\quad F=\tilde F_1\cdots\tilde F_N
\end{displaymath}
with 
\begin{displaymath}
\tl W_k={\rm diag}\{I_{\phi_k},W_k,I_{\eta_k}\},\; k=1,\dots,N-1;\;
\tilde F_k={\rm diag}\{I_{\phi_k},F_k,I_{\chi_k}\},\; k=1,\dots,N,
\end{displaymath}
for suitable sizes $\phi_k$, $\eta_k$ and $\chi_k$. Here $W_k\;(k=1,\dots,N-1)$ are $(n_k+\rho_{k-1})\times(n_k+\rho_{k-1})$
and $F_k\;(k=1,\dots,N)$ are $(s_{k-1}+m_k)\times(s_{k-1}+m_k)$ unitary
matrices.
\end{remark}

\bibliographystyle{amsplain}

\end{document}